\documentclass[12pt]{amsart}
\usepackage{amssymb}
\usepackage{euscript}
\let\script\EuScript
\let\cal\mathcal

\newtheorem{theorem}{Theorem}
\newtheorem{lemma}[theorem]{Lemma}
\newtheorem{corollary}[theorem]{Corollary}
\newtheorem{proposition}[theorem]{Proposition}
\newtheorem{remark}[theorem]{Remark}

\newtheorem{definition}[theorem]{Definition}

\def\eqref#1{(\ref{eq#1})}
\def\eqlabel#1{\label{eq#1}}

\def\D{\script D}

\def\N{\mathbb N}

\def\R{\mathbb R}
\def\M{\cal{M}}
\let\<\langle

\def\T{\tau}

\numberwithin{equation}{section}
\numberwithin{theorem}{section}

\let\\\cr
\let\epsilon\varepsilon
\let\phi\varphi

\def\supp{\operatorname{supp}}

\def\max{\operatorname{max}}

\title[Haagerup $L^p$-spaces]{Kadec-Pe\l czy\'nski
Decomposition
for Haagerup $L^{\boldsymbol p}$-spaces}

\author{Narcisse Randrianantoanina}
\address{Department of Mathematics and Statistics, Miami University, Oxford,
Ohio 45056}
\thanks{Supported in part by NSF Grant DMS-9703789}
\email{randrin@muohio.edu}

\begin{document}

\begin{abstract}
Let $\M$ be a von Neumann algebra (not necessarily semi-finite).  We provide
a generalization of the classical Kadec-Pe\l czynski 
subsequence decomposition of
bounded sequences in $L^p[0,1]$ to the case of the Haagerup $L^p$-spaces $(1\le
p<\infty)$.  In particular, we prove that if $(\varphi_n)_n$
is a bounded sequence in the predual
$\M_*$  of $\M$, then there exist  a subsequence 
$(\varphi_{n_k})_k$ of $(\phi_n)_n$, a decomposition
$\varphi_{n_k}= y_k+ z_k$ such that $\{y_k,\ k\ge 1\}$ is 
relatively weakly
compact and the support projections $s(z_k)\downarrow_k 0$ 
(or similarly mutually disjoint).  
As an application, we prove that
every non-reflexive subspace of the  dual of any given 
$C^*$-algebra (or Jordan triples)
 contains asymptotically isometric
copies of $\ell^1$ and therefore fails the fixed
point property for nonexpansive mappings. These generalize
earlier  results for the case of preduals of 
semi-finite von Neumann algebras.
\end{abstract}

\maketitle

\section{Introduction}
In \cite{KP}, Kadec and Pe\l czy\'nski proved 
a fundamental property that
if $1\leq p<\infty$ then every bounded sequence
$(f_n)$ in $L^p[0,1]$ has a
subsequence that can be decomposed  into
two extreme sequences $(g_k)$ and $(h_k)$, where the $h_k$'s are pairwise
disjoint,  the $g_k$'s are $L_p$-equi-integrable that is
$\lim\limits_{m(A) \to 0}\sup\limits_{k}\Vert \chi_{A} g_{k}
 \Vert_p \rightarrow 0$ and $h_k \perp g_k$ for every 
$k \geq 1$.
This result is generally known as the Kadec-Pe\l czy\'nski
 subsequence decomposition and has been investigated by several authors for the cases of
Banach lattices and symmetric spaces (see for instance
\cite{JMST}  and \cite{WE}).

Motivated by the characterization of relatively weakly compact
subsets of preduals of von Neumann algebras by Akemann
\cite{AK}, the above decomposition was studied in
\cite{D3LRS} for
non-commutative $L^1$-spaces  associated with
  semi-finite
von Neumann algebras equipped with  distinguished, faithful,
normal, semi-finite traces. A more general situation  on
$E(\M,\T)$, where $E$ is a symmetric space of functions on
$(0, \infty)$ and $\M$ is a semi-finite
von Neumann algebra, was studied in \cite{Ran10}. In particular,
the result in \cite{D3LRS} was generalized for $L^p(\M,\T)$ 
for all $0 <p <\infty$.

The aim of the present paper is  to provide extensions of the
Kadec-Pe\l czy\'nski decomposition theorem for general
von Neumann algebras which are not necessarily semi-finite. There are 
many different methods of constructions of non-commutative 
$L^p$-spaces associated with  arbitrary 
von Neumann algebras; for instance,
those of Araki-Masuda \cite{AM}, Haagerup \cite{HAA},
Hilsum \cite{HIL}, Izumi \cite{IZ}, Kosaki \cite{KOS}, Terp
\cite{TER2} and many others. But it is known that, for a given
von Neumann algebra $\M$ and a fixed index $p$, all these 
$L^p$-spaces are isometrically isomorphic. 
We will consider  Haagerup's $L^p$-spaces since they can 
be viewed as spaces of operators that can be embedded as
subspaces of symmetric spaces of measurable operators
obtained from  semi-finite von Neumann algebras via 
crossed  product (see a 
brief description below). 
Our main result is Theorem~4.1 which roughly says that
any bounded sequence in $L^p(\M)$ has a subsequence that
can be splitted into two sequences; one is uniformly integrable
and the other consists of elements supported by decreasing 
projections that converges
to zero.
Our initial motivation is the case $p=1$ 
where $L^1(\M)$ can simply be
viewed as the predual of $\M$. This case allows us to get
informations on  
 copies of $\ell^1$ in duals of $C^*$-algebras.
It has been known that every non-reflexive subspace
of duals of $C^*$-algebras contains complemented copies of 
$\ell^1$ \cite{PF2}. 
On the other hand, Dowling and Lennard 
showed in \cite{DL1}  that for 
$L^1[0,1]$, these complemented copies can be chosen 
to be asymptotically isometric. 

Using the main  decomposition for the case $p=1$,
we can conclude that every non reflexive
subspace of duals of $C^*$-algebras contains
 sequences that generate complemented copies of 
$\ell^1$ and  are  
 asymptotically isometric.
As in \cite{D3LRS} and \cite{DL1}, these asymptotically 
isometric copies of $\ell^1$ yield self maps on
convex bounded sets that fail to have any fixed points. 
These lead to a more general structural
consequence  that non-reflexive subspaces
of duals of $JB^*$-triples fail the fixed point property for
self-maps on closed bounded convex sets.

The paper is organized as follows: 
in Section~2 below, we set some preliminary background on
Haagerup $L^p$-spaces. In particular, we provide a brief 
discussion on its connection to the semi-finite case and 
define the notion of uniformly integrable sets
in these $L^p$-spaces. Section~3
is devoted to the proof a key result which is essentially
the crusial part of the paper. We present our main results
in Section~4 and finally, Section~5 is where we provide all 
the applications on copies of $\ell^1$ on duals of
$C^*$ algebras and $JB^*$-triples.

Our notation and terminology are standard as may be found in 
\cite{D1} for Banach spaces, \cite{KR} and \cite{TAK} for
operator algebras.

\section{ Non-commutative $L^p$-spaces }

In this section, we will describe  different spaces involved
and discuss some properties that will be crusial for 
the presentation.  We will begin from the semi-finite case.
We denote by
${\cal N}$ a semi-finite von Neumann algebra  on a 
Hilbert space $\script H$, with a
distinguished normal,
faithful semi-finite trace $\tau$.  The identity in
$\cal N$ will be denoted
by $\bf{1}$. 
 A closed and densely defined operator $a$ on $\script{H}$ is said to be
affiliated with $\cal N$ if $ua = au$ for all unitary  operator $u$ in the
commutant $\cal{N}'$ of $\cal N$.

A closed and densely defined operator $x$, affiliated with 
$\cal N$,
 is called
$\tau$-measurable if for every $\epsilon > 0$, there exists an orthogonal
projection $p \in \cal{N}$ such that $p(\script{H})\subseteq
\text{dom}(x)$, $\tau ({\bf 1} -p)<\epsilon$ and $xp \in \cal{M}$.
The set of all $\tau$-measurable operators will be denoted by
$\widetilde{\cal N}$.  The set $\widetilde{\cal N}$ is a
${*}$-algebra with respect to the strong
sum, the strong product and the adjoint operation.
 Given a self-adjoint
operator $x$ in $\widetilde{\cal N}$ and $B$ a Borel
subset of $\R$,  we denote by $\chi_{B}(x)$ the projection
$\int_{B} 1\ de^x$ where 
$e^{x}(\cdot)$ is the  spectral
measure of $x$.  
For fixed $x \in \widetilde{\cal N}$ and $t \geq 0$,
we recall
$$\mu_{t}(x) = \inf \left\{s\geq 0 : \tau(e^{\vert x \vert}(s,
\infty))\leq t \right\}.$$
  The function  $\mu_{(.)}(x):[0, \infty) \rightarrow [0,
\infty]$ is called the generalized singular value function (or decreasing
rearrangement) of $x$. 
For a complete study of $\mu_{(.)}$, we refer 
to \cite{FK}.

If $E$ is a 
 symmetric  (r.i.  for short) quasi-Banach
function space on $\R^+$,  the symmetric space of measurable
operators $E(\cal{N}, \tau)$  is defined by setting
$$E(\cal{N},\tau) := \left\{x \in \widetilde{\cal M} : \mu(x)\in E 
\right\}$$
and
$$\left\Vert x \right\Vert_{E(\cal{N}, \tau)} =
 \left\Vert \mu(x) \right\Vert_{E} \  \ \text{for \ all} \ x \in
E(\cal{N},\tau).
$$
The space $E(\cal{N}, \tau)$ is
a (quasi) Banach space and  is
often referred to as the non-commutative version
of the (quasi) Banach function space $E$.  We remark that if
$0<p<\infty$ and $E = L^p
(\R^{+}, m)$  then
$E(\cal{N}, \tau)$ coincides with the usual
non-commutative $L^p$-space associated to the semi-finite von Neumann
algebra $\cal N$.  
 We refer to \cite{DDP1}, \cite{DDP3} and \cite{X} for
extensive background on the space $E(\cal{N},\tau)$.

\medskip

We now provide a short description of the Haagerup $L^p$-spaces.
Let assume that $\M$ is a general von Neumann algebra (not necessarily
semi-finite). 
Let $\cal N$  be the crossed product of $\M$ by the modular
automorphism group $(\sigma_t)_{t \in \R}$ of a fixed 
semi-finite weight on $\M$.
 The von Neumann algebra
$\cal N$  admits the dual action $(\theta_s)_{s \in \R}$ and a
normal faithful semi-finite trace $\T$ satisfying,
$\T\circ \theta_s =e^{-s}\T$, $s \in \R$. 
For $1\leq p< \infty$, the Haagerup $L^p$-space associated with 
$\M$ is defined by
 $$L^p(\M) := \{ x \in \widetilde{\cal{N}}: \theta_s(x)=e^{-s/p}x, s \in
 \R\}.$$

It is well known that there is a linear order isomorphism
$\phi \to h_\phi$ from $\M_*$ onto $L^1(\M)$. One can define
a positive linear functional $Tr$ on $L^1(\M)$ by setting
$$Tr(h_\phi)=\phi({\bf 1}),\ \quad \phi \in \M_*.$$
For $1\leq p < \infty$, 
the spaces $L^p(\M)$ are Banach spaces with the norm 
 defined by
$$ \Vert x \Vert_p = \left(Tr(|x|^p)\right)^{\frac1p},
\ \quad  \text{for}\ x\in L^p(\M).$$
 For complete details on  the construction of $L^p(\M)$, we refer to
 \cite{TE}. Also, it was shown in \cite[Lemma~4.8]{FK}  that if
 $x \in L^p(\M)$, $1\leq p <\infty$, then
  $$\mu_t(x)=t^{-1/p}||x||_p, \ \ \  t>0.$$
where the singular value is relative to the canonical trace on
$\cal{N}$.

We recall that if $1\leq p < \infty$, then the Lorentz space
$L^{p,\infty}(\R^{+},m)$ is the set of (class of) all Lebesgue
measurable functions
on $\R^{+}$ with the norm
 $$ ||f||_{p,\infty}=\sup_{t>0}\{t^{1/p}\mu_t(f)\}.$$
 It is well known that if $1< p < \infty$, then the space
 $L^{p,\infty}(\R^+,m)$ equipped with the  equivalent
 Calderon norm given by
 $$||f||_{(p,\infty)}=\sup_{t>0}\left\{t^{1/p-1}\int_0^t \mu_s(f)\
  ds \right\}, \ \ \ f \in  L^{p,\infty}(\R^+,m),$$
is a symmetric Banach function space on $\R^+$ 
with the Fatou property. The following proposition is 
an immediate  
consequence of the above remarks.
\begin{proposition}
If $1<p<\infty$, then the space $L^p(\M)$ is a closed subspace
of the symmetric space of measurable operators
$L^{p,\infty}(\cal{N},\T)$. Moreover if $1/q +1/p=1$, then
 $$ ||x||_p = q||x||_{(p,\infty)}$$
for all $x \in L^p(\M)$.
\end{proposition}

\medskip

Let us now extend the notion of uniform integrability to 
the Haagerup $L^p$-spaces.
Following \cite{D3LRS}, we define uniform integrability
in $L^p(\M)$ as in Akemann's characterization of
relatively weakly compact subsets of $\M_*$.

\begin{definition}
Let  $1\leq p <\infty$
and $K$ be a bounded subset of $L^p(\M)$.
We  say that $K$ is
uniformly integrable 
 if $\lim\limits_{n \to \infty} \sup\limits_{\phi \in K}
\left\Vert e_n \phi e_n \right\Vert = 0$
for every decreasing sequence $(e_n)_n$ of projections 
in $\M$ with
$e_n \downarrow_n  0$.
\end{definition}

We note that for $p=1$, a subset $K$ is uniformly 
integrable $L^1(\M)$
if and only if it is relatively weakly compact.

\medskip
Throughout, $\script{D}$ denotes the set of all sequences of decreasing
projections in $\M$ that converges to zero;
$\script{D}: = \left\{(e_n)_{n};\  \text{the $e_n$'s are projections in $\M$ and}
 \  e_n \downarrow_n 0 \right\}$. Also for any subset 
$K$ of $L^p(\M)$, $|K|$ denotes the set of all modudi of elements of
$K$; $|K|:=\left\{|x|; x \in K\right\}$. 

\medskip

\noindent
{\bf Fact~1.}  {\it If $x\in L^p(\M)$ and $y\in L^p(\M)$ are such that $x\perp y$ (i.e.
($\supp x)\perp (\supp y)$) then $||x+y||^p=||x||^p+||y||^p$.}

\begin{proof}  $||x||^p=Tr(|x|^p)$ and if $x\perp y$ as
elements of $\cal N$, $|x+y|^p =|x|^p+|y|^p$ and
therefore  $||x+y||^p =Tr (|x|^p+|y|^p)=
||x||^p+||y||^p$.
\end{proof}

\medskip
\noindent
{\bf Fact~2.} {\it If $x \in L^p(\M)$ and $e$ is  a projection
in $\M$, then $\Vert x\Vert^p \geq \Vert exe \Vert^p +
\Vert (1-e)x(1-e) \Vert^p$.}
\begin{proof}
Set $u=2e-1$. It is clear that $u \in \M$ is unitary and
$exe + (1-e)x(1-e) =\frac{1}{2}\left(x +uxu^*\right)$. 
It follows that $\Vert exe + (1-e)x(1-e) \Vert^p \leq 
\Vert x \Vert^p$ and hence
$\Vert exe \Vert^p + \Vert (1-e)x(1-e) \Vert^p \leq 
\Vert x \Vert^p$.
\end{proof} 

We finish this section with 
the following two lemmas which can be proved using 
similar arguments as in the semi-finite case 
(\cite{D3LRS}, \cite{Ran10}) and
will be used in the sequel. Details are left to the readers.

\begin{lemma}  Let  $1\leq p<\infty$,
 $(p_n)_n \in \script{D}$ 
and $K$ be a bounded subset of $L^p(\cal{M})$ such that
for each $n_0 \geq 1$, the sets
$({\bf 1} - p_{n_0})K$ and $\vert K({\bf 1}-p_{n_0}) \vert$ 
are uniformly integrable.
Then $K$ is uniformly integrable if and only if
$\lim\limits_{n \to \infty}\sup_{\phi \in K} 
\Vert p_n\phi p_n \Vert =0.$
\end{lemma}

\begin{lemma} Let $1\leq p<\infty$,
$(\phi_n)_n$ be a bounded sequence in   $L^p(\M)$ 
and $(p_n)_n \in \script{D}$.
Assume that $\lim\limits_{n \to \infty} 
\sup\limits_{k} \left\Vert p_n \phi_k p_n
\right\Vert=\gamma > 0$
 then there exists a subsequence $(\varphi_{k_n})$ so
that $\lim\limits_{n \to \infty} 
\left\Vert p_n\phi_{k_n}p_n \right\Vert
 = \gamma$.
\end{lemma}

\section{Preliminary Results}
This section is devoted to the proof of Theorem~3.1 below
which is the key result that we will use to prove our main
theorem. We remark that the case of finite von-Neumann
algebras can be obtained with minor changes from the proof
of the commutaive case (see \cite{D3LRS}).

\begin{theorem}
Let $\M$ be a $\sigma$-finite von-Neumann algebra and $1 \leq p < \infty$.
Assume that $K$ is a  subset of 
the positive part of the unit ball of $L^p(\M)$ 
that is not uniformly integrable. 
 Then there exists a sequence 
$(\varphi_n)_n \subset K$ and $(f_n)_n
\in \D$ such that:
$$
\sup\left\{\lim_{n\to\infty}
\sup_{k \in \N} \Vert e_n \varphi_k e_n\Vert;\
(e_n)_n \in \D \right\} = \lim_{n \to \infty}
\sup_{k \in \N} \Vert f_n \varphi_k f_n \Vert > 0. $$
\end{theorem}

\begin{lemma}  Let $\cal N$ be a semi-finite von Neumann algebra with
distinguished faithful normal semi-finite trace $\T$ as above 
and $E$
be a symmetric 
 quasi-Banach function space on $(0,\infty)$.
If $x\in E(\cal N,\T)$ and $u\in\cal N$ then
$$
||xu||_{E}\le ||x||^{\frac 12}_{E}
 \cdot ||u^*|x|u||^{\frac 12}_{E} \le
||x||^{\frac 34}_{E}\cdot ||uu^*|x|uu^*||^{\frac
14}_{E}.$$
\end{lemma}

\begin{proof}  Let $x=v|x|$ be the polar decomposition of $x$.  Then
\begin{equation*} 
\begin{split}
\left\Vert xu \right\Vert =||v|x|u|| &\le|||x|u||
=|||x|^{\frac 12}|x|^{\frac 12}u|| \\
&\le |||x|^{\frac 12}||_{E^{(2)}}\cdot 
|||x|^{\frac 12}u||_{E^{(2)}}\\
&=\Vert x \Vert^{\frac 12}_{E} \cdot \Vert u^*|x|u 
\Vert^{\frac 12}_{E}.
\end{split}
\end{equation*}
For the second inequality,
\begin{equation*}
\begin{split}
\left\Vert xu \right\Vert &\leq
 ||x||^{\frac 12}_{E}\cdot 
|||x|^{\frac 12}u||_{E^{(2)}}  \\
&=||x||^{\frac 12}_{E}\cdot 
||u^*|x|^{\frac 12}||_{E^{(2)}} \\
&\le ||x||^{\frac 12}_E\cdot |||x|^{\frac 12}uu^*|x|^{\frac 12}||^{\frac
12}_{E}                                            \\
&\le||x||^{\frac 12}_E 
\left(\left\Vert|x|^{\frac 12}uu^*
\right\Vert_{E^{(2)}}\cdot \left\Vert|x|^{\frac
12}\right\Vert_{E^{(2)}}\right)^{\frac 14}                               \\
&=||x||_E^{\frac 12}\cdot ||uu^*|x|uu^*||^{\frac 14}_E\cdot ||x||^{\frac
14}_E                                           \\
&=||x||^{\frac 34}_E \cdot ||uu^*|x|uu^*||^{\frac 14}_E.
\end{split}
\end{equation*}
\end{proof}

 Lemma~3.2 shows in particular that if $1\leq p<\infty$, $x\in L^p(\M)$ 
and $u\in \M$ then
$$||xu||\le||x||^{\frac 34}\cdot 
||uu^*|x|uu^*||^{\frac 14}.$$

\begin{lemma}
Let $\gamma >0$ and $(\phi_k)_k$ be a sequence in the positive part
of the unit ball of $L^p(\M)$.  If there exists a sequence $(a_n)_n$ in the
unit ball of $\M$ with $a_n\downarrow_n 0$ and such that
$\lim_{n\to
\infty}\sup_k ||a_n\phi_ka_n||\ge\gamma$.  Then for every 
$\varepsilon>0$, there exists a
sequence $(s_n)_n$ of projections with:
\begin{itemize}
\item[(i)] $s_n \leq s_1$ for every $n \geq 1$;
\item[(ii)] $s_n \to 0$ for the strong operator topology;
\item[(iii)] for every $n_0 \in \N$, $\lim_{n\to\infty}\sup_k
\Vert(s_{n_0}a_n s_{n_0})\phi_k(s_{n_0}a_n s_{n_0})\Vert
\ge\gamma-\varepsilon$.
\end{itemize}
\end{lemma}

\begin{proof}
  Fix $\delta>0$ with $\delta \leq (\epsilon/8)^2$
and define the sequence of projections as follows:
\begin{equation*}
\begin{cases}
s_1:=\chi_{(\delta,1)}(a_1) &and \\
 s_n:=\chi_{(\delta,1)}
(s_1a_n s_1) &{for\  n\ge 2}.
\end{cases}
\end{equation*}
 Clearly  $s_n$ is a subprojection of the
support of $s_1 a_n s_1$ so  $s_n \leq s_1$.
Also
 $\delta s_n\le s_n(s_1a_n s_1)s_n$, and since
$s_n$ and $s_1a_n s_1$ are commuting operators,
$\delta s_n\le s_n(s_1a_n s_1)s_n \le s_1 a_n s_1$ and therefore
$s_n \to 0$ so $(i)$ and $(ii)$ are verified. 

\medskip
\noindent
{\it Claim: Let  $n_0 \in \N$ and $n\geq n_0$, for every
$\phi \in L^p(\M)$ with $\Vert \phi \Vert \leq 1$,
$\Vert \phi a_n(1-s_{n_0})\Vert \leq 2\delta^{1/2}$. 
Similarily, 
$\Vert \phi(1-s_{n_0})a_n\Vert \leq 2\delta^{1/2}$.}

\medskip

To see this claim, it is enough to notice that 
\begin{equation*}
\begin{split} 
\Vert \phi a_n (1-s_{n_0}) \Vert &\leq
\Vert \phi a_n s_1(1-s_{n_0}) \Vert +
\Vert \phi a_n (1-s_1) \Vert \\
&\leq
\Vert \phi \Vert \cdot \Vert a_n s_1(1-s_{n_0}) \Vert_\infty +
\Vert \phi \Vert \cdot \Vert a_n (1-s_1) \Vert_\infty \\
&\leq \Vert(1-s_{n_0}) s_1 a_{n}^2 s_1(1-s_{n_0})
 \Vert_{\infty}^{1/2} +
\Vert  (1-s_1)a_{n}^2 (1-s_1) \Vert_{\infty}^{1/2}\\
&\leq \Vert(1-s_{n_0}) s_1 a_{n} s_1(1-s_{n_0})
 \Vert_{\infty}^{1/2} +
\Vert  (1-s_1)a_{n} (1-s_1) \Vert_{\infty}^{1/2}
\end{split}
\end{equation*}
and since $(a_n)_n$ is a decreasing sequence and 
$1\leq n_0 \leq n$, we get
$$\Vert \phi a_n (1-s_{n_0}) \Vert \leq
\Vert(1-s_{n_0}) s_1 a_{n_0} s_1(1-s_{n_0})
 \Vert_{\infty}^{1/2} +
\Vert  (1-s_1)a_1 (1-s_1) \Vert_{\infty}^{1/2}
\leq  2\delta^{1/2}.$$
A similar estimate can be established for 
$\Vert \phi(1-s_{n_0})a_n\Vert$ which verifies the claim. 

To complete the proof, let $\phi \in L^p(\M)$, 
$\Vert \phi \Vert \leq 1$. For $n \geq n_0$, we can write
$a_n\phi a_n$ as: 
$$
a_n\varphi a_n = (s_{n_0}a_ns_{n_0}) \varphi a_n 
+s_{n_0} a_n (1-s_{n_0})\varphi a_n
+ (1-s_{n_0})a_n \varphi a_n $$
and using the claim above,
$\Vert a_n \phi a_n \Vert \leq
\Vert (s_{n_0}a_ns_{n_0}) \varphi a_n \Vert + 4\delta^{1/2}$.
A similar estimate would give
$\Vert (s_{n_0}a_n s_{n_0})\phi a_n \Vert \leq
\Vert (s_{n_0}a_ns_{n_0}) \varphi (s_{n_0}a_n s_{n_0}) \Vert 
+ 4\delta^{1/2}$ and 
combining these two estimates, we get
$$\Vert a_n \phi a_n \Vert \leq
\Vert (s_{n_0}a_ns_{n_0}) \varphi (s_{n_0}a_n s_{n_0}) \Vert 
+ 8\delta^{1/2}.$$
This shows that
$$\lim_{n \to \infty}\sup_k
\Vert (s_{n_0}a_ns_{n_0}) \varphi_k (s_{n_0}a_n s_{n_0}) \Vert 
\geq \gamma -8\delta^{1/2} \geq \gamma-\epsilon.$$
 The proof is complete.
\end{proof}

\medskip
The next result shows that using projections in the definition of 
uniform integrability is not essential. One can use
elements of the positive part of the unit ball of $\M$.

\begin{proposition}  Let $\gamma >0$ and $(\phi_k)_k$ be a sequence in the positive part
of the unit ball of $L^p(\M)$.  If there exists a sequence $(a_n)_n$ in the
unit ball of $\M$ with $a_n\downarrow_n 0$ 
and such that
$\lim_{n\to
\infty}\sup_k ||a_n\phi_ka_n||\ge\gamma$.  Then for every $\varepsilon>0$, there exists a
sequence $(p_n)_n\in\D$ with 
$p_n \leq supp(a_1)$ for all $n \geq 1$ and 
such that $\lim_{n\to\infty}
\sup_{k}
||p_n\phi_kp_n||\ge\gamma-\varepsilon$.
\end{proposition}

\begin{proof} The sequence $(p_n)$ will be constructed 
inductively. Let $(\epsilon_j)_j$ be a sequence in the 
open interval $(0,\epsilon)$ such that 
$\sum_{j=1}^\infty \epsilon_j =\epsilon$ and $\omega_0$ be a 
faithful state in $\M_*$.

By Lemma~3.3, one can choose a sequence of projections
$(s_{n}^{(1)})_n$ with $ s_{n}^{(1)} \leq s_{1}^{(1)}$ 
for every $n \geq 1$,  $s_{n}^{(1)} \to 0$ (as n tends to 
$\infty$)
satisfying the conclusion of Lemma~3.3 for $(a_n)_n$,
$\gamma$ and $\epsilon_1$.

Choose $n_1 \geq 1$ such that  
$\omega_0\left(s_{n_1}^{(1)}\right)\leq 1/2$.
From (iii) of Lemma~3.3,
$$\lim_{n\to\infty}\sup_k
\Vert(s_{n_1}^{(1)}a_n s_{n_1}^{1})\phi_k
(s_{n_1}^{(1)}a_n s_{n_1}^{(1)})\Vert
\ge\gamma-\varepsilon_1.$$
Reapplying  Lemma~3.3,  on 
$(a_{n}^{(2)})_n=(s_{n_1}^{(1)}a_n s_{n_1}^{(1)})_n$,
$\gamma-\epsilon_1$ and $\epsilon_2$, one would get a sequence
of projections 
$(s_{n}^{(2)})_n$ with $ s_{n}^{(2)} \leq s_{n_1}^{(1)}$ 
for every $n \geq 1$,  $s_{n}^{(2)} \to 0$ 
(as n tends to infinity).
As above, on can choose $n_2$ such that 
$\omega_0\left(s_{n_2}^{(2)}\right) \leq 1/{2^2}$ and
$$\lim_{n\to\infty}\sup_k
\Vert(s_{n_2}^{(2)}a_n s_{n_2}^{(2)})\phi_k
(s_{n_2}^{(2)}a_n s_{n_2}^{(2)})\Vert
\ge\gamma-\epsilon_1 -\epsilon_2.$$
The induction is clear, repeating the argument above would give
a decreasing sequence of projections
$s_{n_1}^{(1)} \ge s_{n_2}^{(2)} \ge \cdots \ge
s_{n_j}^{(j)} \ge \cdots$ so that  for every
$j \geq 1$, $\omega_0\left(s_{n_j}^{(j)}\right)\leq 1/{2^j}$
and
$$\lim_{n\to\infty}\sup_k
\Vert(s_{n_j}^{(j)}a_n s_{n_j}^{(j)})\phi_k
(s_{n_j}^{(j)}a_n s_{n_j}^{(j)})\Vert
\ge\gamma-\sum_{i=1}^j\epsilon_i.$$
If for every $j\ge 1$, we set $p_j= s_{n_j}^{(j)}$ then 
$(p_j)_j$ belongs to $\D$ and 
$$\sup_k
\Vert(p_ja_j p_j)\phi_k
(p_ja_j p_j)\Vert
\ge\gamma-\sum_{i=1}^j\epsilon_i $$
which shows that
$$\lim_{j\to \infty} \sup_k \Vert(p_ja_j p_j)\phi_k
(p_ja_j p_j)\Vert
\ge\gamma-\epsilon $$
and since $\Vert p_ja_j \Vert_\infty \leq 1$, the
desired conclusion follows.
\end{proof}

\begin{proposition}
Let $K$ be as in the statement of 
Theorem~3.1.  There exists a sequence 
$(\varphi_k)_k$ in
$K$ such that 
$$\sup\left\{\lim_{n\to\infty}
\sup_{k \in \N} \Vert e_n \varphi_k
e_n\Vert;  \
(e_n)_n \in \D\right\} 
= \sup\left\{\varliminf_{n\to\infty} \Vert 
e_n \varphi_n e_n\Vert;\ 
(e_n)_n \in \D \right\} > 0.$$
\end{proposition}

\begin{proof}  Set
$\alpha_0 := \sup\left\{
\lim_{n\to\infty}
\sup_{\phi \in K} \Vert e_n \phi e_n\Vert;\
(e_n) \in \D \right\}$
and let $(\varepsilon_j)_j $ be a subset of the open interval
 $(0,1)$ such that
$\Pi^\infty_{j=1}(1 -\varepsilon_j) >0$.

Since $\alpha_0 > 0$, one can choose a sequence
$(y_n)_n$ in $K$ and $(e_n^{(1)})_n \in \D$ such that
$$\lim_{n\to\infty}
\sup_{k \in \N} \Vert e_n^{(1)}y_k
e_n^{(1)} \Vert \geq \alpha_0 (1 - \varepsilon_1).$$ 
 A further
subsequence $(y_k^{(1)})_k \subset (y_k)$ can be chosen so that
$$\lim_{n\to\infty}\Vert
e_n^{(1)}y_n^{(1)} e_n^{(1)} \Vert 
\geq \alpha_0(1 - \varepsilon_1).$$
Set $\displaystyle{\alpha_1:= \sup \left\{\lim_{n\to\infty}
\sup_{k \in \N} \Vert e_ny_k^{(1)}
e_n\Vert ; (e_n) \in \D \right\}}$.  It is clear that $\alpha_1 \geq
\alpha_0 (1 - \varepsilon_1)$ and as above a sequence $(y_k^{(2)})_{k \geq
1} \subseteq (y_k^{(1)})_k$ can be chosen so that
$$\lim_{n\to\infty}\Vert
e_n^{(2)}y_n^{(2)} e_n^{(2)} \Vert \geq \alpha_1(1 - 
\varepsilon_2).$$
Inductively, one can construct 
sequences $(y_n)_n \supseteq (y_n^{(1)})_n \supseteq
(y_n^{(2)})_n \supseteq \ldots (y_n^{(j)})_n \supseteq 
\ldots$  in $K$ and sequences
$(e_n^{(1)})_n, (e_n^{(2)})_n,\ldots, (e_n^{(j)})_n ,\ldots$ in $\D$ so that
for every $j \geq 1$,
 $$\lim_{n\to\infty}\Vert
e_n^{(j)}y_n^{(j)} e_n^{(j)} \Vert \geq 
\alpha_{j-1}(1 - \varepsilon_j).$$
Let $(\phi_n)_n$ be the diagonal sequence obtained from $(y_n^{(j)})_n, j \geq
1$.  For every $j \geq 1$, $(\phi_n)_{n \geq j}$ is a subsequence of
$(y_n^{(j)})_{n \geq 1}$ so 
$$\displaystyle{\varliminf_{n\to\infty}\Vert
e_n^{(j)}\phi_n e_n^{(j)} \Vert \geq \alpha_{j-1}(1 - 
\varepsilon_j)}$$
 and
$$\displaystyle{\sup\left\{\lim_{n\to\infty}
         \sup_{k\in\N}
         \Vert e_n \phi_k e_n \Vert; (e_n) \in \D\right\} 
\leq \alpha_j}.$$
We note that $\alpha_{j-1} \geq \alpha_j \geq \alpha_{j-1} (1 -
\varepsilon_j)$ so for every $j \geq 1$,
\begin{equation*}
\begin{split}
\sup\left\{\lim_{n\to\infty}\sup_{k\in\N}
\Vert e_n \phi_k e_n \Vert ; (e_n) \in \D\right\} 
&\leq \alpha_j \cr
&\leq \alpha_j (1 - \varepsilon_{j+1}) \frac1{1-
         \varepsilon_{j+1}}\cr
&\leq \frac 1 {1 - \varepsilon_{j+1}} \varliminf_{n \to \infty}
         \Vert e_n^{(j+1)} \phi_n e_n^{(j+1)} \Vert
\end{split}
\end{equation*}
which implies that
$$
\sup\{\lim_{n\to\infty}\sup_{k\in\N} \Vert e_n \phi_k
e_n \Vert;\
(e_n)_n\in \D\}
\le\frac 1{1-\varepsilon_{j+1}}
\sup \left\{\varliminf_{n\to \infty}\Vert e_n\phi_ne_n||;\
(e_n)_n\in  \D \right\}.$$

Taking the limit as $j$ goes to $\infty$,
$$
\sup\left\{\lim_{n\to\infty}\sup_{k\in  \N} \Vert e_n \phi_ke_n
\Vert;\ (e_n)_n\in \D \right\}
\le\sup\left\{\varliminf_{n\to\infty}
\Vert e_n\phi_ne_n \Vert;\ (e_n)_n\in \D \right\}.$$
The other inequality is trivial.

To check that
$\sup\left\{\varliminf_{n\to\infty}
\Vert e_n\phi_ne_n \Vert;\ (e_n)_n\in \D \right\}>0$, it
is plain  that
$$
\varliminf_{n\to\infty}\Vert e_n^{(j)}\phi_ne_n^{(j)}
\Vert \ge\alpha_{j-1}(1-\epsilon_j)
\ge\alpha_1\Pi_{j=2}^{\infty}(1-\epsilon_j)>0.$$

The proof of the proposition is complete.
\end{proof}

\medskip
\noindent
{\bf Proof of Theorem~3.1}

Let $(\phi_n)$ be the sequence in $K$ obtained from 
Proposition~3.5
;  i.e. $(\phi_n)_n$ is the sequence in $K$ satisfying:
$$
\sup\left\{\varliminf_{n\to\infty} 
\Vert e_n\phi_n e_n \Vert;\ (e_n)_n\in  \D\right\}
=\sup\left\{\lim_{n\to\infty}\sup_{k\in \N}
\Vert e_n\phi_ke_n \Vert;\
(e_n)_n\in  \D \right\}:=\alpha > 0
$$

\noindent
{\it Claim:  $\alpha$ is attained.}

\medskip

Assume the opposite i.e. for every $(p_n)_n\in \D$,
$\varliminf_{n\to\infty}\Vert p_n \phi_np_n \Vert<\alpha$.

Inductively, we will construct sequences of integers and
projections in $\M$ satisfying the following conditions:
\begin{equation}\eqlabel{1}
m_1\le m_2\le\cdots\le m_j\le\cdots\  \text{a sequence in 
 $\N$};
\end{equation}

\begin{equation}\eqlabel{2}
n_1<n_2<\cdots<n_j<\cdots\  \text{infinite sequence in $\N$};
\end{equation}

sequences $(p_n^{(1)})_{n\ge 1}, (p_n^{(2)})_{n\ge 1},\cdots$ 
in $\D$
such that for every $j\ge 2$ and every $n\ge 2$,
\begin{equation}\eqlabel{3}
p_n^{(j)}\perp \sum_{k=1}^{j-1}p_{n_k}^{(k)};
\end{equation}
if we set $ (f_n^{(1)})_{n\ge 1}=(p_n^{(1)})_{n\ge 1}$ and
\begin{equation}\eqlabel{4}
f_n^{(j)}=\cases f_n^{(j-1)} \vee p_{n_{j-1}}^{(j)}\quad n< n_{j-1}\\
f_n^{(j-1)}+p_n^{(j)}\qquad n\ge n_{j-1}\endcases
\end{equation}
then
\begin{equation}\eqlabel{5}
\lim_{n\to\infty}\sup_{k\in \N}
\Vert f_n^{(j)}\phi_k f_n^{(j)} \Vert\ge\alpha
(1-\frac 1{2^{m_j-1}}),
\end{equation}

\begin{equation}\eqlabel{6}
\sup_{k\in \N}
\Vert f_{n_j}^{(j)}\phi_kf_{n_j}^{(j)} \Vert<\alpha
(1-\frac 1{2^{m_j}})
\end{equation}
and
\begin{equation}\eqlabel{7}
\varliminf_{n\to\infty}\Vert f_n^{(j)}\phi_nf_n^{(j)}
\Vert^p \ge
\varliminf_{n\to\infty}\Vert f_n^{(j-1)}\phi_nf_n^{(j-1)}
\Vert^p+
\frac{\alpha^{4p}}{(2^{m_{j-1}+3})^{4p}}.
\end{equation}

\medskip

Fix a sequence  $(p_n^{(1)})_{n\ge 1} \in \D$ such that
$\varliminf_{n\to\infty}
\Vert p_n^{(1)}\phi_n p_n^{(1)}\Vert\ge\alpha(1-\frac
1{2^2})$ and choose $m_1\in \N$ such that
$$
\alpha(1-\frac 1{2^{m_1-1}})\le\lim_{n\to\infty}\sup_{k\in
\N} \Vert p_n^{(1)}\phi_kp_n^{(1)}\Vert
<\alpha(1-\frac 1{2^{m_1}})
$$
(such $m_1$ exists since $\alpha$ is not attained).  

Assume that the 
construction is done for $1, 2, \cdots, (j-1)$.  
By the definition of
$\alpha$, one can choose $(q_n)_n\in\script D$ so that
$\varliminf_{n\to\infty}\Vert q_n \phi_n q_n \Vert>
\alpha(1-\frac{1}{2^{m_{j-1}+1}})$. Writing 
$q_n\phi_n q_n$ in the form
$$
q_n \phi_nq_n = q_nf^{(j-1)}_{n_{j-1}}\phi_nf^{(j-1)}_{n_{j-1}}q_n
  +q_nf^{(j-1)}_{n_{j-1}}\phi_n(1-f^{(j-1)}_{n_{j-1}})q_n +
 q_n(1-f_{n_{j-1}}^{(j-1)})\phi_n q_n,
$$
one can see that
$$
\Vert q_n \phi_nq_n\Vert \le
\Vert f_{n_{j-1}}^{(j-1)}\phi_nf_{n_{j-1}}^{(j-1)}\Vert +
2\Vert \phi_n(1-f_{n_{j-1}}^{(j-1)})q_n\Vert.
$$
Applying Lemma~3.2  for $x=\phi_n$ and 
$u=(1-f_{n_{j-1}}^{(j-1)})q_n$, we get
\begin{equation*}
\begin{split}
\Vert q_n\phi_nq_n \Vert
&\le \Vert f_{n_{j-1}}^{(j-1)}\phi_nf_{n_{j-1}}^{(j-1)}
\Vert +\\
&
2\Vert \phi_n\Vert^{\frac34} \cdot
\Vert(1-f_{n_{j-1}}^{(j-1)})q_n(1-f_{n_{j-1}}^{(j-1)})
\phi_n(1-f_{n_{j-1}}^{(j-1)})
q_n(1-f_{n_{j-1}}^{(j-1)})\Vert^{\frac 14}.
\end{split}
\end{equation*}
 Applying \eqref{6} for $(j-1)$ gives
$$
\Vert q_n\phi_nq_n \Vert
\le\alpha (1-\frac 1{2^{m_{j-1}}})
+ 2\Vert(1-f_{n_{j-1}}^{(j-1)})q_n(1-f_{n_{j-1}}^{(j-1)})\phi_n
(1-f_{n_{j-1}}^{(j-1)})q_n(1-f_{n_{j-1}}^{(j-1)})
\Vert^{\frac14}.
$$
Taking the limit (as $n$ tends to $\infty$),
$$
\alpha (1-\frac 1{2^{m_{j-1}+1}})\le\alpha (1-\frac 1{2^{m_{j-1}}})
+2\varliminf_{n \to \infty}
\Vert(1-f_{n_{j-1}}^{(j-1)})q_n(1-f_{n_{j-1}}^{(j-1)})\phi_n
(1-f_{n_{j-1}}^{(j-1)})q_n(1-f_{n_{j-1}}^{(j-1)})
\Vert^{\frac14}.$$
which implies that
$$\varliminf_{n\to\infty}
\Vert(1-f_{n_{j-1}}^{(j-1)})q_n(1-f_{n_{j-1}}^{(j-1)})
\phi_n(1-f_{n_{j-1}}^{(j-1)})q_n(1-f_{n_{j-1}}^{(j-1)})
\Vert \ge\frac{\alpha^4}{(2^{m_{j-1}+2})^4}.
$$
If we set $a_n^{(j)}=(1-f_{n_{j-1}}^{(j-1)})q_n(1-f_{n_{j-1}}^{(j-1)})$
then $a_n^{(j)}\downarrow_n 0$ and
$$
\varliminf_{n\to\infty}\Vert a_n^{(j)}\phi_na_n^{(j)}\Vert\ge
\frac{\alpha^4}{(2^{m_{j-1}+2})^4}.
$$
Applying Proposition~3.4 for
$(\phi_n)_n$, $(a_n^{(j)})_n$,
$\gamma=\frac{\alpha^4}{(2^{m_{j-1}+2})^4}$  and  
$\epsilon=
\frac{\alpha^4}{(2^{m_{j-1}+2})^4} - 
\frac{\alpha^4}{(2^{m_{j-1}+3})^4}$ 
would provide a sequence $(p_n^{(j)})\in \D$ such that
$$
\varliminf_{n\to\infty}\Vert p_n^{(j)}\phi_np_n^{(j)}\Vert\ge
\frac{\alpha^4}
{(2^{m_{j-1}+3})^4}.
$$
Since $p_n^{(j)}\le \supp(a_1^{(j)})\le 
1-f_{n_{j-1}}^{(j-1)}$,
it is clear that $p_n^{(j)} \perp f^{(j-1)}_{n_{j-1}}$ for every $n \geq
1$ so \eqref{3} is verified.

If we define $(f^{(j)}_n)$ as in \eqref{4} then
appropriate $m_{j} \geq m_{j-1}$ and $n_{j} > n_{j-1}$ can be choosen so
that \eqref{5} and \eqref{6} would be satisfied.

Now since $p^{(j)}_n + f^{(j-1)}_n = f_n^{(j)}$ for $n \geq n_j$,

\begin{equation*}
\begin{split}
\Vert f_n^{(j)} \phi_n f_n^{(j)} \Vert^p
&\geq \Vert f_n^{(j-1)} \phi_n f_n^{(j-1)} \Vert^p + 
 \Vert p_n^{(j)} \phi_n p_n^{(j)} \Vert^p \cr
&\geq \Vert f_n^{(j-1)} \phi_n f_n^{(j-1)} \Vert^p
+ \frac {\alpha^{4p}} {(2^{m_{j-1}+3})^{4p}}
\end{split}
\end{equation*}
and \eqref{7} is verified. The construction is done.

\medskip

To complete the proof of the theorem, we note from 
\eqref{7} that

\begin{equation*}
\begin{split}
\varliminf_{n \to \infty} \Vert f^{(j)}_n \phi_n 
f^{(j)}_n \Vert^p \geq 
\varliminf_{n \to \infty} \Vert f^{(1)}_n 
\phi_n f^{(1)}_n \Vert^p + \alpha^{4p} \sum^{j-1}_{k=1}
\frac {1} {(2^{m_{k}+3})^{4p}}
\end{split}
\end{equation*}

So the series $\displaystyle\sum^{\infty}_{k=1}
\frac {1} {(2^{m_{k}+3})^{4p}}$ is
convergent. In particular
$\lim_{k \to \infty} m_k = \infty$.

We note from \eqref{4} that every $j \geq 1$ and $n \geq n_{j-1}$,
$f^{(j)}_n =\displaystyle\sum^j_{k=1} p^{(k)}_n =
 \vee^j_{k=1} p_n^{(k)}$;
this is the case since $(n_{j})$ is increasing so if $n \geq n_{j-1}$ then
$n \geq n_{l-1}$ for all $l \leq j$ and hence
\begin{equation*}
\begin{split}
f^{(j)}_n &= f^{(j-1)}_n + p^{(j)}_n \cr
&= f^{(j-2)}_n + p^{(j-1)}_n + p^{(j)}_n \cr
&= \sum^j_{k=1}p^{(k)}_n
\end{split}
\end{equation*}
and all the $(p_n^{(k)})_{1 \leq k \leq j}$ are 
mutually disjoint.  

Now
choose an increasing sequence $(k_{j})$ so that $k_j > 
\max(k_{j-1},n_{j-1})$,  
$\omega_0 (f^{(j)}_{k_{j}}) < \frac {1}{2^j}$ and
\begin{equation}\eqlabel{8}
\alpha (1-\frac {1}{2^{m_{j}{-2}}}) \leq
\sup_{k \in \N} \Vert 
f^{(j)}_{k_{j}}  \phi_k  f^{(j)}_{k_{j}} \Vert
\end{equation}
(this last condition is possible from \eqref{5}).

\medskip
\noindent
{\it Claim:
$\{f^{(j)}_{k_{j}};j \geq 1\}$ 
is a commuting family of projections in $\cal M$.}

\medskip
In fact for each $j \geq l$,
$f^{(j)}_{k_{j}} = \sum^j_{k=1} p^{(k)}_{k_{j}}$ and 
$f^{(l)}_{k_{l}} = \sum^l_{k=1} p^{(k)}_{k_{l}}$.
For each $1 \leq k\leq l$, 
$p^{(k)}_{k_{l}} \geq p^{(k)}_{k_{j}}$ and
$p^{(k)}_{k_{l}} \perp \sum\limits_{s=1;s\neq k}^j
p^{(s)}_{k_{j}}$, hence

$$f^{(j)}_{k_{j}} f^{(l)}_{k_{l}} = f^{(l)}_{k_{l}}
 f^{(j)}_{k_{j}} =
\sum^l_{k=1} p^{(k)}_{k_{j}}$$
and the claim follows.
\medskip

Set $\cal{S}$ to be a maximal abelian von Nemann subalgebra
of $\M$ generated by
$\left\{f^{(j)}_{k_{j}}\ ; j \geq 1 \right\}$.  
Since $\cal{S}$ is abelian, $\omega_0$ restricted to 
$\cal{S}$ is a faithful tracial state on $\cal{S}$.

Set $p_n: =\displaystyle\operatornamewithlimits \vee_{j \geq n}
f^{(j)}_{k_{j}}$ (where the supremum is taken in $\cal{S}$).
It is clear that
$$\omega_0 (p_{n}) \leq \sum^{\infty}_{j=n} \omega_0
(f^{(j)}_{k_{j}}) \leq \sum^{\infty}_{j=n} \frac {1}{2^{j}}$$

so $\omega_0(p_n)\to 0$ which shows that $(p_n)_n 
\in \D$.
Moreover, since $p_n \geq  f^{(n)}_{k_{n}}$
for all $n\geq 1$,  condition~\eqref{8} implies that
$$ \sup_{k \in \N} \Vert p_{n}
\phi_k p_n \Vert  \geq 
\alpha (1- \frac {1}{2^{m_n-2}}).$$ 
This would show  that
$\lim_{n \to \infty}
\sup_{k \in \N} \Vert p_{n} \phi_{k}
p_{n} \Vert = \alpha$.

This is a contradiction with the initial assumption that 
$\alpha$
is not attained. The proof is complete.

\qed

\section{Main Result}

The main results in this section are Theorem~4.1 and
Theorem~4.4 which generalize the classical 
Kadec-Pe\l czy\'nski subsequence decomposition to 
bounded sequences in the Haagerup $L^p$-spaces.

\begin{theorem}
Let $\M$ be a von Nemann algebra, $1 \leq p < \infty$ and $(\varphi_n)_n$
be a bounded sequence in $L^p(\M)$. Then there exist a subsequence
$(\varphi_{n_{k}})_{k}$ of $(\varphi_{n})_{n}$, 
two bounded sequences
$(y_{k})$, $(z_{k})$ in $L^{p}(\M)$ and a decreasing 
sequence of projections
$e_k \downarrow_k 0$  in $\M$ such that:

\begin{itemize}
\item [(i)] $\varphi_{n_{k}} = y_k + z_k$ for all $k \geq 1$;
\item [(ii)] $\{y_k, \ k \geq 1 \}$ is uniformly integrable 
in $L^p(\M)$;
\item [(iii)] $z_k = e_k z_k  e_k$ for all $k \geq 1$.
\end{itemize}
\end{theorem}

\begin{proof}
We will assume first that $\M$ is $\sigma$-finite.  Without loss of
generality, we can and do assume that $\Vert \varphi_n \Vert \leq 1$
for all
$n \geq 1$ and $\{\varphi_n, n \geq 1 \}$ is not uniformly integrable.  We
will show that there exist a sequence $(n_k)$ in $\N$ and $(e_k)_k \in \D$
such that 
the bounded set $\{\varphi_{n_{k}} - e_k \ \varphi_{n_{k}}\ e_k ; k \geq 1\}$
is uniformly integrable in $L^p(\M)$.

By Theorem~3.1, there exists a subsequence
of  $(\varphi_{n})_{n}$ ( which we will denote again by
$(\phi_n)_n$) and
$(e_{n})_n \in \D$ such that

\begin{equation*}\begin{split}
\sup &\{\displaystyle\operatornamewithlimits \lim_{n \to \infty}
\displaystyle\operatornamewithlimits \sup_{k \in \N}
\Vert f_n (\vert \varphi_k \vert + \vert \varphi^{*}_{k} \vert) f_n \Vert ;
(f_n) \in \D\} \cr
&=\displaystyle\operatornamewithlimits 
\lim_{n \to \infty}
\displaystyle\operatornamewithlimits \sup_{k \in \N}
\Vert e_n (\vert \varphi_k \vert + 
\vert \varphi^{*}_{k} \vert) e_n \Vert =
\alpha >0.
\end{split}\end{equation*}

Choose a subsequence $(\varphi_{n_{k}})_{k} \subset (\varphi_{k})$ so that

\begin{equation}\eqlabel{1*}
\lim_{k \to \infty} \Vert e_{k} \
\varphi_{n_{k}} \ e_{k} \Vert = \alpha.
\end{equation}

Set $u_{k}: = \varphi_{n_{k}}$ and $v_{k}: = 
u_k - 
e_{k}  u_k  e_{k}$
for all $k\geq 1$.  

\medskip

\noindent
{\it Claim:  The set $V = \{v_k ; \ k \in \N \}$ is uniformly
integrable in $L^{p}(\M)$.}

\medskip

To see this claim, we will first prove the following 
intermediate lemma:
\begin{lemma}
Let $n_{0} \in \N$, $(1- e_{n_{0}}) V$ and $\vert V (1-e_{n_{0}}) \vert$ are
uniformly integrable subsets of $L^p(\M)$.
\end{lemma}

We will show that  $\vert V (1-e_{n_{0}}) \vert$ is uniformly integrable.
Assume the opposite.  There exists $(f_n)_{n} \in \D$ such that
$$ \lim_{n \to \infty}
 \sup_{k \in \N} \Vert f_n \vert
v_k (1- e_{n_{0}}) \vert f_n \Vert > 0. $$
From this, 
there  would exists $(p_{n})_{n} \in \D$ with
$p_n \leq 1 - e_{n_{0}}$  and such that
$$\lim_{n \to \infty}\sup_{k \in \N} \left\Vert p_n (\vert
u_k \vert + \vert u^{*}_{k} \vert p_{n}) \right\Vert > 0.$$

In fact, for each $k \geq 1$, 
if we denote by $\omega_k$  the partial isometry in $\M$
so that
$\vert v_k (1-e_{n_{0}}) \vert = \omega_k v_k (1-e_{n_{0}})$,
then
\begin{equation*}
\begin{split}
\Vert f_n \vert v_k (1- e_{n_{0}}) \vert f_n \Vert =
 \Vert f_n  \omega_k 
v_k  (1-e_{n_{0}}) f_n \Vert 
\leq \Vert v_k (1-e_{n_{0}}) f_n \Vert
\end{split}
\end{equation*}

Note that for $k \geq n_0$, $e_k(1-e_{n_{0}})=0$ so 
$\Vert f_n \vert v_k
(1-e_{n_{0}}) \vert f_n \Vert \leq \Vert u_k (1-e_{n_{0}}) 
f_n \Vert$
and by Lemma~3.2,
$$\Vert f_n \vert v_k (1-e_{n_{0}}) \vert f_n \Vert \leq \Vert u_k
\Vert^{\frac{3}{4}} \cdot \Vert (1-e_{n_{0}}) f_n (1-e_{n_{0}}) \vert u_k
\vert (1-e_{n_{0}}) f_n (1-e_{n_{0}})\Vert^{\frac{1}{4}}$$
which shows that
$$\Vert f_n \vert v_k (1-e_{n_{0}}) \vert f_n \Vert \leq 
\Vert (1-e_{n_{0}}) f_n (1-e_{n_{0}}) \vert u_k \vert (1-e_{n_{0}}) f_n
(1-e_{n_{0}}) \Vert^{\frac {1}{4}}.$$

Let $a_n = (1-e_{n_{0}}) f_n (1-e_{n_{0}})$.
It is clear that $ a_n \downarrow_n 0$ and using
Proposition~3.4, we conclude that there exists 
$(p_n)_{n} \in \D$ , $p_n \leq
1-e_{n_{0}}$ such that 
$\lim_{n \to\infty}\sup_{k \in \N} \Vert p_n
\vert u_k \vert p_n \Vert > 0$. In particular:
$$
\lim_{n \to \infty}
\sup_{k \in \N}
\Vert p_n (\vert u_k \vert + \vert u^{*}_{k} \vert) 
p_n \Vert > 0.$$

Now choose a subsequence $(k_{j}) \subseteq \N$ so that
 there exists $\gamma>0$ satisfying

\begin{equation}\eqlabel{2*}
\lim_{j \to \infty}
\Vert p_j (\vert u_{k_{j}} \vert + \vert u^{*}_{k_{j}} 
\vert) p_j \Vert = \gamma > 0.
\end{equation}

Since $p_j \leq 1-e_{n_{0}}$ for all $j$, 
$e_{k_{j}} \perp p_j$ for $k_j >
n_0$ and therefore
\begin{equation*}\begin{split}
\left\Vert (e_{k_{j}} + p_j) \left(\vert u_{k_{j}} \vert + 
\vert u^{*}_{k_{j}}\vert \right)
(e_{k_{j}} + p_j) \right\Vert^{p} 
\geq \left\Vert e_{k_{j}} \left(\vert u_{k_{j}} 
\vert + \vert u^{*}_{k_{j}} \vert \right)
e_{k_{j}} \right\Vert^{p} 
+ \left\Vert p_j \left(\vert u_{k_{j}} \vert + \vert u^{*}_{k_{j}} 
\vert \right) p_j \right\Vert^{p}
\end{split}\end{equation*}
and taking the limit as $j \to \infty$, \eqref{1*} and
\eqref{2*}
imply  $\alpha^{p} \geq \gamma^{p} +
\alpha^{p}$. This is a contradiction   since $\gamma > 0$.
The proof of the lemma is complete.

\medskip

To complete the proof of the theorem, assume that $V$ is not uniformly
integrable.  Using Lemma~2.2  and Lemma~4.2,
$$\lim_{n \to \infty}\sup_{k \in \N}
\Vert e_{n} v_{k} e_n \Vert > 0.$$
Again, choose a subsequence $(k_{n}) \subseteq \N$ so that
\begin{equation}\eqlabel{3*}
\lim_{n \to \infty}
\Vert e_n v_{k_{n}} e_n \Vert > 0.
\end{equation}

\medskip
{\it Claim:
$\left\Vert e_n v_{k_{n}} e_n \right\Vert^{2} \leq 4 
\left\Vert (e_n - e_{k_{n}})
\left(\vert u_{k_{n}} \vert + 
\vert u^{*}_{k_{n}} \vert\right) (e_n - e_{k_{n}})
\right\Vert$.}

\medskip
To see this claim, we note that since $e_n \geq e_{k_{n}}$,
$e_n v_{k_{n}} e_n = e_n u_{k_{n}} e_n - 
e_{k_{n}} u_{k_{n}} e_{k_{n}}$  so
$e_n v_{k_{n}} e_n = (e_n - e_{k_{n}}) u_{k_{n}} e_n + e_{k_{n}}
u_{k_{n}} (e_n - e_{k_{n}})$ and therefore
\begin{equation*}\begin{split}
\Vert e_n v_{k_{n}} e_n \Vert &\leq \Vert (e_n - e_{k_{n}}) u_{k_{n}} \Vert +
\Vert u_{k_{n}} (e_n-e_{k_{n}}) \Vert \cr
&\leq \Vert u^{*}_{k_{n}} (e_n -e_{k_{n}}) \Vert + \Vert u_{k_{n}} (e_n
- e_{k_{n}}) \Vert \cr
&\leq \Vert u^{*}_{k_{n}} \Vert^{\frac{1}{2}} \cdot \Vert (e_n-e_{k_{n}}) \vert
u^{*}_{k_{n}} \vert (e_n -e_{k_{n}}) \Vert^{\frac{1}{2}} 
+ \Vert u_{k_{n}} \Vert^{\frac{1}{2}} \cdot \Vert (e_n - e_{k_{n}}) \vert
u_{k_{n}} \vert (e_n - e_{k_{n}}) \Vert^{\frac{1}{2}}
\end{split}\end{equation*}

and since $\Vert u_{k_{n}} \Vert \leq 1$,
\begin{equation*}
\begin{split}
\Vert e_n v_{k_{n}} e_n \Vert
&\leq \Vert (e_n - e_{k_{n}})
\vert u^{*}_{k_{n}} \vert (e_n - e_{k_{n}}) 
\Vert^{\frac{1}{2}} 
+  \Vert (e_n - e_{k_{n}}) \vert u_{k_{n}} \vert (e_n - e_{k_{n}})
\Vert^{\frac{1}{2}} \cr
&\leq  2 \Vert (e_n - e_{k_{n}}) \left(\vert u_{k_{n}} \vert +
\vert u^{*}_{k_{n}} \vert \right) (e_n - e_{k_{n}}) \Vert^{\frac{1}{2}}
\end{split}
\end{equation*}
and the claim follows.

\medskip

From the claim above and equation \eqref{3*}, there exists
$\nu >0$ such that

\begin{equation}\eqlabel{4*}
\varliminf_{n \to \infty} \left\Vert (e_n - e_{k_{n}}) 
\left(\vert u_{k_{n}} \vert + \vert
u^{*}_{k_{n}} \vert \right) (e_n - e_{k_{n}}) \right\Vert = \nu > 0.
\end{equation}

 Observe that since $e_{k_{n}} \perp (e_n -
e_{k_{n}})$,
\begin{equation*}
\begin{split}
\Vert e_n \left(\vert u_{k_{n}} \vert + 
\vert u^{*}_{k_{n}}\vert\right) e_n \Vert^{p} 
\geq \Vert e_{k_{n}} \left(\vert u_{k_{n}} \vert + 
\vert u^{*}_{k_{n}} \vert\right)
e_{k_{n}} \Vert^{p} 
+ \Vert (e_n - e_{k_{n}}) \left(\vert u_{k_{n}} \vert + \vert
u^{*}_{k_{n}} \vert\right) (e_n - e_{k_{n}}) \Vert^{p}.
\end{split}
\end{equation*}
Taking the limit ( as $n \to \infty$) together with
\eqref{1*} and \eqref{4*} would imply
 $\alpha^{p} \geq  \nu^{p} +\alpha^{p}$.
This is  a contradiction
since $\nu >0$.

By setting $y_k := v_k$ and $z_k := e_k u_k e_k$,  
the proof for the $\sigma$-finite case is complete.

\bigskip

For the general case, let $\M$ be a von Neumann algebra (not necessarily
$\sigma$-finite) and $(\varphi_{n})_{n}$ in $L^{p}(\M)$ 
as in the theorem. 
Fix an  orthogonal
family of cyclic projections  $ (e_{\alpha})_{\alpha \in I}$
in $\M$ such that
$1 = \vee_{\alpha \in I} e_{\alpha}$  ( see for instance,
\cite[Proposition 5.5.9, p. 336]{KR1} ) .

\begin{lemma}
There exists a countably decomposable projection $e \in \M$ such that for
all $n \geq 1$, $e \varphi_n = \varphi_{n}e = \varphi_n$.
\end{lemma}

\medskip

\noindent
For each $n \in \N$ and $\epsilon > 0$, set
$E_{n, \epsilon}: = \{\alpha \in I \ ; \ \Vert e_{\alpha}\varphi_n \Vert >
\epsilon \}$ and
$E_n : = \{\alpha \in I \ ; \ \Vert e_{\alpha} \varphi_{n}
 \Vert \neq 0 \}$.

\medskip

\noindent
{\it Claim:
$E_{n,\epsilon}$ is finite (hence $E_n$ is countable). }

\medskip

To see this,   assume that
$E_{n,\epsilon}$ is infinite.  Then there exists an infinite sequence
$(e_k)_{k} \subset (e_{\alpha})_{\alpha \in I}$ such that $\Vert e_k
\varphi_n \Vert > \epsilon$ for all $k \in \N$.  If $J$ is a finite subset
of $\N$, then
\begin{equation*}
\begin{split}
\Vert \sum_{k \in J} e_k \varphi_n \Vert &= \Vert
(\sum_{k \in J} e_{k}) \varphi_n \Vert \cr
&= \Vert(\vee_{k \in J} e_k) \varphi_n
\Vert \leq \Vert \varphi_n \Vert.
\end{split}
\end{equation*}
So $\Vert\sum_{k \in J} e_k \varphi_n \Vert \leq \Vert
\varphi_n \Vert$  (a constant independant of $J$)
which shows that  $\sum_{k
=1}^{\infty}e_k \varphi_n$ is a weakly unconditionally 
Cauchy (w.u.c.) series in
$L^p (\M)$ but since $L^p(\M)$ does not contain any copies of
$c_0$, $\sum_{k
=1}^{\infty}e_k \varphi_n$ is  unconditionally 
convergent and hence 
$\lim_{k \to \infty} \Vert e_k
\varphi_n \Vert = 0$ 
(see for instance  \cite{D1} p.45). This is in  contradiction with 
the assumption $ \Vert e_k \varphi_n \Vert \geq
\epsilon$ for all $k \in \N$. We proved that
$E_{n,\epsilon}$ is finite. It is clear that $E_n =
\cup_{k \in \N} E_{n,
\frac{1}{k}}$ so it is at most countable. The claim is
verified.

\medskip

Similarly, if $R_n = \{\alpha \in I , \Vert \varphi_n e_{\alpha} \Vert
\neq 0 \}$ then $R_n$ is at most countable.

Let $C = \displaystyle\operatornamewithlimits\cup^{\infty}_{n = 1}(R_n \cup
E_n)$ ; $C$ is at most countable and if
$e =  \vee_{\alpha \in C} e_{\alpha}$
then $e$ is the union of a countable family of disjoint  
cyclic projections in $\M$ so
$e$ is countably decomposable in $\M$ 
(\cite[Proposition~ 5.5.19  p.340 ]{KR1}).  The
construction of $e$ implies that $e \varphi_n = \varphi_n e = \varphi_n$
for all $n \geq 1$. The lemma is proved.

To conclude the proof of the theorem, 
consider the von Neumann algebra $e \M e$.  Since
$e$ is countably decomposable, $e \M e$ is $\sigma$-finite.  Let $T: \M \to
e \M e$ be the map that takes $x \in \M$ to $exe$. 
The map  $T$ is bounded and is
 weak$^{*}$ to weak$^{*}$ continuous so there exists a
map $S: (e \M e)_{*} \to \M_{*}$ so that $S^{*} =T$.

Let $R:\M_{*} \to (e \M e)_{*}$ be the restriction map.  The operators $T$ and
$R$ can be interpolated and since $L^{p}(\M)$ (resp. 
$L^{p}(e \M e)$) is
isometrically isomorphic to $(\M, \M_{*})_{\theta}$ (resp. $(e \M e,(e \M
e)_{*})_{\theta})$ for $\theta = \frac{1}{p}$, 
(see \cite{TER2}), we get a bounded linear map 
$T_p:
L^p (\M) \to L^p(e \M e)$. 
Similarly, if one considers the inclusion map $e \M e
\to \M$ and $S:(e \M e)_{*} \to \M_{*}$ as above, 
then by interpolation, we obtain a map 
$S_p:
L^{p}(e \M e) \to L^{p}(\M)$.

Apply the $\sigma$-finite case to the
sequence $(T_p(\varphi_n))_{n \geq 1}$ in 
$ L^{p} (e \M e)$ 
to get 
 a decomposition
$$T_{p} (\varphi_{n_{k}}) = y_k + z_k  \quad  \forall \ k 
\geq1$$
with $(y_k)_k$ and $(z_k)_k$ satisfying the conclusion of 
the theorem.
It is enough to consider the decomposition:
$$\varphi _{n_{k}} = S_p (y_k) + S_p(z_k) 
\quad \forall \ k \geq 1.$$

The proof  is complete.
\end{proof}

\medskip
The theorem which follows shows that, as in the semi-finite 
case, the decreasing projections in Theorem~4.1 can be replaced
by mutually orthogonal projections. Its proof is identical 
to that of the semi-finite case (\cite{Ran10}, Theorem~3.7).

\begin{theorem}
Let $\M$ be a von Neumann algebra and $1\leq p <\infty$, 
  Let $(\phi_n)_{n}$ be a bounded sequence in $L^p(\cal{M})$ then there exists a
subsequence $(\phi_{n_{k}})$ of $\phi_n$, 
bounded sequences $(y_k)$ and
$(\zeta_k)_{k}$ in $L^p(\cal{M})$ and mutually orthogonal
 sequence of projections
$(e_k)_{k}$ in $\M$ such that:
\begin{itemize}
\item [(i)] $\phi_{n_{k}} = y_k + \zeta_k$  for all $k \geq 1$;
\item [(ii)] $\{y_k:\  k\geq 1\}$ is uniformly integrable and
$e_k y_k e_k = 0 $ for all $\ k \geq 1$;
\item [(iii)] $(\zeta_k)_{k}$ is such that $e_k \zeta_k e_k = \zeta_k$
 for all $k \geq1$.
\end{itemize}
\end{theorem}

\begin{remark}
For $1<p<\infty$, it should be noted that since $L^p(\M)$ is
is a closed subspace of $L^{p,\infty}(\cal{N},\T)$ and
$L^{p,\infty}(\R^+,m)$ has the Fatou property, one could apply
the semi-finite case of the the Kadec-Pe\l czy\'nski subsequence
decomposition to any bounded sequence of $L^p(\M)$ 
(viewed as bounded 
sequence in $L^{p,\infty}(\cal{N},\T))$. However, that
procedure would provide decreasing projections in $\cal{N}$ and
as is noted in \cite[Remarks~3.5 (iii)]{Ran10}, these 
projections are either of finite trace or their orthogonal
complements are of finite trace which guaranties that
projections obtained from applying the semifinite case 
cannot be in $\M$.
\end{remark}

\section{Applications}

A  result of Maurey (\cite{Mau}, see also \cite{JPOS})
states that every reflexive subspace of $L^1[0,1]$ has the fixed
point property for nonexpansive mappings (FPP). Later,
Dowling and Lennard showed that the converse of 
Maurey's result is valid: every non-reflexive 
subspace $L^1[0,1]$ fails the FPP
 (\cite{DL1}).
This section  is for the study of 
 generalizations to the case of 
duals of $C^*$-algebras and  requires the notion of asymptotically 
isometric copies of $\ell^1$ which was introduced by 
Dowling and Lennard in \cite{DL1}.

\begin{definition}
A Banach space $X$ is said to contain asymptotically isometric
copies of $\ell^1$ if for every null sequence $(\epsilon_n)$
of positive numbers, there exists a sequence $(x_n)$ in $X$
such that:
$$
\sum_{n=1}^\infty (1-\epsilon_n)|a_n|
\leq \left\Vert \sum_{n=1}^\infty a_nx_n \right\Vert
\leq  \sum_{n=1}^\infty |a_n|.
$$
for all $(a_n) \in \ell^1$.
\end{definition}

The following result is a generalization of \cite{D3LRS}.
\begin{theorem}
let $\cal{A}$ be a $C^*$-algebra. Every non-reflexive
subspace of $\cal{A}^*$ contains asymptotically
isometric copies of $\ell^1$.
\end{theorem}
\begin{proof} Note that $\cal{A}^{**}$ is a von Neumann
algebra so subspaces of $\cal{A}^*$ are subspaces of 
preduals of von Neumann algebras. The proof then follows 
the argument used in \cite{D3LRS} using Theorem~4.4.
Details are left to the readers.
\end{proof}

\begin{remark}
In \cite{BEL}, B\'elanger proved an improved version of 
the Akemann's characterization of weak compactness on
preduals of von Neumann algebras. He then went on to 
show that non-reflexive preduals of von Neumann algebras
contain complemented copies of $\ell^1$.
This fact can also be deduced from a result of Pfitzner
\cite{PF2} which states that $C^*$-algebras have 
Pe\l czy\'nski property~(V)
so their duals have property~(V*) (\cite{PL1}). It is plain from
Theorem~4.4 that the asymptotically isometric copies of
$\ell^1$ in Theorem~5.2 are complemented with
good projection constants.
\end{remark}

For the next extension, we recall that 
 $JB^*$-triples  are all those Banach spaces whose 
 open unit balls
are bounded symmetric domains \cite{UP}. Examples of 
$JB^*$-triples are
$C^*$-algebras and Hilbert spaces. Other important
examples are the so-called {\it Cartan factors}
$C^k (k=1,2,\dots,6)$  where the rectangular Cartan factor
$C^1=\cal{L}(H,K)$ consists of bounded operators between 
Hilbert spaces, the symplectic factor $C_{n}^2$ is 
$\{z \in \cal{L}(H);\ z=-jz^*j\}$ where $j:H \to H$ is a conjugate
linear isometric involution, the Hermitian Cartan
factor $C^3$ is $\{z \in \cal{L}(H);\ z=jz^*j\}$,
$C^4$ is the spin factor, $C^5$ is the (finite dimentional)
exceptional Cartan factor consisting of $1\times 2$ matrices
over the complex Caley numbers $\bf O$ and
$C^6$ is the set of all $3\times 3$ Hermitian matrices over 
$\bf O$. Dual $JB^*$-triples are called $JBW^*$-triples.
For more informations, we refer to \cite{CI2}, 
\cite{HOR1} and \cite{HOR2}.

\begin{corollary}
If $\cal J$ is a $JB^*$-triple then every non-reflexive 
subspace of $\cal{J}^*$ contains asymptotically isometric
copies of $\ell^1$.
\end{corollary}

For the poof we will need two  lemmas on stability of asymptotically
isometric copies of $\ell^1$.
\begin{lemma}
Let $E_1$ and $E_2$ be weakly sequentially complete
Banach spaces so that any sequence
equivalent to the unit vector
basis of $\ell^1$ in $E_j$ $(j=1,2)$ has a normalized
block that is asymptotically isometric to $\ell^1$ then
every sequence  equivalent to the unit vector basis of 
$\ell^1$ in $E_1 \oplus_1 E_2$ has a normalized block
that is asymptotically isometric to $\ell^1$.
\end{lemma}
\begin{proof} Let $\{U_n=(x_n , y_n)\}_{n=1}^\infty$
be a sequence in $E_1 \oplus_1 E_2$ that is 
equivalent to $\ell^1$. After taking subsequences,
either $(x_n)_n$ or $(y_n)_n$ is equivalent to $\ell^1$.
Let assume that $(x_n)_n$ is equivalent to $\ell^1$.
We have two cases.

\noindent
{\it Case~1: The sequence 
$(y_n)_n$ is weakly convergent.} By taking normalized 
blocks, 
we can assume that $(x_n)_n$ is asymptotically isometric
to $\ell^1$ and 
$\lim_{n\to \infty} \Vert y_n \Vert =0$.
There exists a null sequence $(\epsilon_n)$
of positive numbers
such that:
$$
\sum_{n=1}^\infty (1-\epsilon_n)|a_n|
\leq \left\Vert \sum_{n=1}^\infty a_nx_n \right\Vert
\leq  \sum_{n=1}^\infty |a_n|.
$$
for all $(a_n) \in \ell^1$ but since
$\Vert \sum_{n=1}^\infty a_n U_n \Vert
=\Vert \sum_{n=1}^\infty a_n x_n \Vert
+\Vert \sum_{n=1}^\infty a_n y_n \Vert$,
we get that 
$$\sum_{n=1}^\infty (1-\epsilon_n) \vert a_n \vert \leq
\left\Vert \sum_{n=1}^\infty a_n U_n \right\Vert
\leq \sum_{n=1}^\infty (1 +
\Vert y_n \Vert)\vert a_n \vert.$$
This concludes that $(U_n)_n $ is 
asymptotically isometric to $\ell^1$.

\noindent
{\it Case~2: The sequence $(y_n)$ is equivalent
to $\ell^1$}.
As above, one can find a block so that both 
the coresponding block for $(x_n)_n$ and $(y_n)_n$
are asymptotically isometric to $\ell^1$.
Set $Z_n:= U_n/2=(x_n/2, y_n/2)$. It can be easily seen
that $(Z_n)_n$ is equivalent to an asympotically 
isometric copy of $\ell^1$ in $E_1 \oplus_1 E_2$.
\end{proof}

\begin{lemma}
Let $(\Omega, \Sigma, \lambda)$ be a measure space and $R$ be a reflexive
Banach space. Every sequence equivalent to the unit vector basis of
$\ell^1$ in $L^1(\lambda, R)$ has a normalized block that is 
asymptotically isometric
to $\ell^1$.
\end{lemma}
\begin{proof}
Let $(f_n)_n$ be a sequence equivalent
 to the $\ell^1$
basis. Since $R$ is reflexive, the sequence $(f_n)_n$
can not be uniformly integrable (see for instance
\cite{DRS}).
Apply the classical Kadec-Pe\l czy\'nski subsequence
decomposition to the sequence 
$(\Vert f_n(\cdot)\Vert)_n$ in $L^1(\lambda)$ to get
a pairwise disjoint sequence 
of measurable sets $(A_n)_n$ such that
$\left\{\Vert f_n(\cdot)\Vert\chi_{\Omega \setminus A_n},
n\geq 1\right\}$ is uniformly integrable.
The space $R$ being reflexive implies that
$\left\{ f_n\chi_{\Omega \setminus A_n},
n\geq 1\right\}$ is relatively weakly compact in 
$L^1(\lambda,R)$. We conclude the proof as in the 
scalar case.
\end{proof} 

\medskip
\noindent{\bf Proof of Corollary~5.4:}
Let $\script{J}$ be a $JB^*$-triple and $X$ be a non-reflexive
subspace of $\script{I}^*$.
Since $\script{J}^{**}$ is a $JBW^*$-triple, we can assume 
that $X$
is a subspace of the predual of a $JBW^*$-triple $\cal{I}$.
By \cite{HOR1} and \cite{HOR2}, $\cal{I}$ admits the following form:
$$
\cal{I}= \left( \sum_\alpha \oplus C(\Omega_\alpha, C^\alpha) 
\right)_{\ell^\infty} \oplus_\infty J^7 \oplus_\infty J^8,
$$
where $C(\Omega_\alpha, C^\alpha)$ is the space of 
continuous functions from a hyperstonean space $\Omega_\alpha$
to a Cartan factor $C^\alpha$,
$J^7=\{a \in M; \Theta(a)=a\} $ with $\Theta:M \to M$ is a 
$w^*$-continuous ${}^*$~-~antiautomorphism of period $2$ on a 
von Neumann algebra $M$ and $J^8$ is a $w^*$-closed
right ideal of a von Neumann algebra $N$.
The predual of $\cal{I}$ is equal to the $\ell^1$-sum
$$
\cal{I}_* = 
 \left( \sum_\alpha \oplus L^1(\Sigma_\alpha, C_*^\alpha) 
\right)_{\ell^1} \oplus_1 J_*^7 
\oplus_1 J_*^8.
$$
By \cite[Theorem~2]{CI2},  the space 
$E_1=\left( \sum_{\alpha\neq 5,6} \oplus L^1(\Sigma_\alpha, C_*^\alpha) 
\right)_{\ell^1} \oplus_1 J_*^7 
\oplus_1 J_*^8 $ is isometric to a
 1-complemented subspace of the predual 
 of a von Neumann algebra so $E_1$ is isometric to 
a subspace of the predual of such von Neumann
algebra and hence satisfies the
assumption of Lemma~5.3. Moreover,
since $C^5$ and $C^6$ are finite dimensional, the space
$E_2 = L^1(\Sigma_5, C^5) \oplus_1 L^1(\Sigma_6, C^6)$ 
satisfies (as does $L^1$-spaces) the assumption  of Lemma~5.3. 
We conclude that 
every sequence  equivalent to the unit vector basis of 
$\ell^1$ in $\cal{I}_*=E_1 \oplus_1 E_2$ has a normalized block
that is asymptotically isometric to $\ell^1$. The proof is 
complete. \qed

\begin{corollary}
If $\cal J$ is a $JB^*$-triple then every non-reflexive 
subspace of $\cal{J}^*$ fails 
the fixed point property
for nonexpansive self-maps on closed
bounded convex sets.
\end{corollary}

\medskip

\noindent
{\bf Acknowledgements.} This project started when
the author was participating in the NSF-supported Workshop
on Linear Analysis and Probability (Summer 1998) at the
Department of Mathematics of the 
Texas A\& M University. The author would like to express
his gratitute to Professor W. Johnson for the invitation and
warm hospitality. The author is also indebted to Professor
P. Dowling for introducing him to the  topic 
of this paper and  for several useful discussions.


\end{document}